\documentclass{amsart}

\usepackage{amsmath}%
\usepackage{amsfonts}%
\usepackage{amssymb}%
\usepackage{graphicx}
\usepackage[all]{xy}
\usepackage{fancybox}
\usepackage{amscd}
\usepackage{latexsym}
\usepackage{wasysym}
\usepackage{pdfpages}
\usepackage{graphicx}

\newcommand{\prova}{\normalfont{\bf Proof: }}
\newtheorem{def1}{Definition}
\newtheorem{teo1}{Theorem}
\newtheorem{exe1}{Example}

\newtheorem{prop1}{Proposition}

\begin{document}
\title[Units of $\mathbb{Z}C_{p^n}$]{Units of $\mathbb{Z}C_{p^n}$}
\author{Raul Antonio Ferraz }
\address[ R. A. Ferraz]
{Instituto de Matem\'atica e Estat\'istica, Universidade de S\~ao Paulo, Brazil}
\email[]{ raul@ime.usp.br}%
\author{ Patr\'icia Massae Kitani }
\address[P. M. Kitani]{ Departamento de Matem\'atica, Universidade Tecnol\'ogica Federal do Paran\'a, Brazil}%
\email[]{kitani@utfpr.edu.br}%
\thanks{$^{\dagger}$ Supported by CAPES and CNPq}
\date{\today}

\begin{abstract}
Let  $p$ be a prime integer and $n,i$ be  positive integers such that \linebreak $S=\{-1, \ \theta, \ \mu_i=1+\theta+\cdots + \theta^{i-1}  \ \mid 1 < i < \frac{p^n}{2} , \ gcd(p^n,i)=1 \}$ generates the group of units of $\mathbb{Z}[\theta],$  where $\theta$ is a primitive ${p^n}$--$th$ root of unity. Denote by $C_{p^n}$ the cyclic group of order $p^n.$ In this paper we describe explicitly a  multiplicatively independent set which generates  a complement to $\pm C_{p^n}$ in  the group of units of the integral group ring of  $C_{p^n}.$
\end{abstract}

\maketitle

\section{Introduction}

Let $\mathbb{Z}G$ be the integral group ring of the group $G$  and let $\mathcal{U}(\mathbb{Z}G)$ be its group of units. 
Many results on the units group of a group rings have been published since it started with Higman in 1940. Due to him, it is well-known that
  if $G$ is finite abelian group,  then $\mathcal{U}(\mathbb{Z}G)=\pm G\times F$ where $F$ is a free abelian group of finite rank. Moreover, the rank is $\frac{|G_0|-2l+m+1}{2},$ where $G_0$ is the torsion subgroup of $G$, $m$ is the number of cyclic subgrupos of $G_0$ of order 2 and $l$ is the number of cyclic subgroups of $G_0.$ (see  \cite{AyoubAyoub},  \cite{Higman} or \cite{Karpilovsky1},).
For a nonabelian finite group it follows from a general result of Borel and Harish-Chandra that the unit group $\mathcal{U}(\mathbb{Z}G)$ is finitely generated. 
The study concerning the group of units of group rings over the integers $\mathbb{Z}G$ is a classic and a difficult problem. Most of the description about $\mathcal{U}(\mathbb{Z}G)$ are given explicitly,  or a subgroup of finite index of $\mathcal{U}(\mathbb{Z}G)$ is given or  a general structure  is studied.

We have provided a subset  of $\mathcal{U} (\mathbb{Z}C_{p^n})$, which is multiplicatively independent and which generates a complement to 
$\pm C_{p^n}$ in $\mathcal{U} (\mathbb{Z}C_{p^n}).$ When $m$ is a general natural number,  few are known about the  set of units  that form  a basis of $\mathcal{U} (\mathbb{Z}C_{m}).$  For $m = 2, 3, 4$ and $6$ we know that  $\mathcal{U} (\mathbb{Z}C_{m})$ is trivial. For $m = 5$ and $m = 8,$ it can be found in \cite{Karpilovsky1} or \cite{Polcino}.
A different description for $m = 8$ was given by Sehgal \cite{Sehgal} using fiber product diagram. Aleev and Panina \cite{Aleev 2} described for $m=7$ and $m=9.$  Ferraz  \cite{Raul}  represented  for prime numbers $m$ between $5$ and $67.$ This last result is very considered in our work. 

Let $\epsilon:\mathbb{Z}G\rightarrow \mathbb{Z}$ be the augmentation homomorphism defined by \linebreak $\epsilon\left( \displaystyle\sum_{g\in G} a_gg\right)=\displaystyle\sum_{g\in G} a_g$ and denote the  subgroup of units with augmentation 1 by $\mathcal{U}_1(\mathbb{Z}G).$ This subgroup of augmentation 1 is called normalized units of $\mathbb{Z}G.$ 

Ferraz \cite{Raul} gave explicitly a multiplicative independent subset of $\mathcal{U}(\mathbb{Z}C_p)$ with generates $\mathcal{U}_1(\mathbb{Z}C_p),$ regarded that $p$ is a prime number such that $$S_{\theta}=\{1+\theta, 1+\theta+\theta^2, \cdots, 1+\theta+\theta^2+\cdots+\theta^{(p-3)/2}\}$$ generates $\mathcal{U}(\mathbb{Z}[\theta]).$ Here, $\theta$ is the primitive $p$--th  root of unity. Let us recall an important unit in \cite{Raul}. Let $C_p=\langle g \rangle,$ $t$ a primitive root  modulo $p$ and  for each $i,$ $1\leq i\leq (p-3)/2,$ then a particular case of the  Hoechsmann unit (see \cite{Sehgal} ) of $\mathbb{Z}C_p$ is 
$$u_i=(\displaystyle \sum_{j=0}^{r-1}g^{tj})(\sum_{j=0}^{t-1}g^{jt^i})-k\widehat{g}=(1+g^t+\cdots +g^{t(r-1)})(1+g^{t^i}+\cdots +g^{t^i(t-1)})-k\widehat{g}$$ where  $r$ is the least positive integer such that $tr$ is congruent to $1$ modulo $p,$  $k$ is the integer $(rt-1)/p$ and $\widehat{g}=1+g+g^2+\cdots +g^{p-1}.$ 
The main result of this paper is the following: 

\begin{teo1}
The set $S=\{u_1, u_2, \cdots, u_{(p-3)/2}\}$ is a multiplicatively independent subset of $\mathcal{U}_1(\mathbb{Z}C_p)$ such that $$\mathcal{U}_1(\mathbb{Z}C_p)=\langle g \rangle \times \langle S \rangle.$$
\end{teo1}

We take the principal idea in the paper of Ferraz and extend to the group ring $\mathbb{Z}C_{p^n}.$ In what follows, we extend the results of section $2$ of \cite{Raul} for powers of  odd prime and after that, in the section 3 to powers of $2.$ In the section 4 we give a set of multiplicatively independent generators of $\mathcal{U}_1(\mathbb{Z}C_{p^n}).$ The last section is to  prove the validity of an hypothesis used in the previous section.
In this paper we give an explicitly set of independent generators of $U(\mathbb{Z}C_{p^n})$ and the fundamental units of our description is the Hoechsmann units.

\section{Odd primes}

In this section we deal with  powers of an odd prime, which we will denote by $p^n$. It is known that there exists an integer $t$ which generates the units of $\mathbb{Z}_{p^n}.$ 
We define $\mu_i=1+\theta+\theta^2+\cdots +\theta^{i-1}$ and  recall the following proposition on Ferraz paper:

\begin{prop1}
Let $p$ be a prime, $5\leq p \leq 67,$ $t\in \mathbb{Z}$ such that $\bar{t}$ generates the group $\mathcal{U}(\mathbb{Z}_p),$  $S_1=\langle-1, \theta, \mu_t, \mu_{t^2}, \cdots, \mu_{t^{(p-3)/2}}\rangle$ and define the following subgroup of $\mathcal{U}(\mathbb{Z}C_p)$ $$S_2=\langle -1, \theta, \mu_t, \mu_t^{-2}\mu_{t^2}, \cdots, \mu_t^{-(p-3)/2}\mu_{t^{(p-3)/2}}\rangle .$$ Then the groups $S_1$ and $S_2$ are equal.    
\end{prop1}

We can extend this result to $p^n.$  

 \begin{prop1}
 Let $p$ be an odd prime, $ I_{p^n}=\{i\in \mathbb{Z}| 1 < i < \frac{p^n}{2} \ \mbox{and} \ gcd(i,p^n)=1\}$ such that $\mathbb{S}=\{-1, \ \theta, \ \mu_i=1+\theta+\cdots +\theta^{i-1} \mid i\in I_{p^n}\}$ generates the group of units of $\mathbb{Z}[\theta]$ and $t\in \mathbb{Z}$ such that $\bar{t}$ generates $\mathcal{U}(\mathbb{Z}_{p^n}).$ Then $$\mathbb{S}_1=\{-1, \theta, \mu_t, \mu_{t^2}, \mu_{t^3}, \cdots, \mu_{t^{\kappa}}\}$$  generates $\mathcal{U}(\mathbb{Z}[\theta])$, where $\kappa=\frac{\phi(p^n)}{2}-1$ and  $\phi$ is the Euler totient. 
 \end{prop1}

\prova

It is easy to see that $\mathbb{S}_1 \subseteq \mathcal{U}(\mathbb{Z}[\theta]).$
Since $\mathcal{U}(\mathbb{Z}[\theta])=\langle \mathbb{S} \rangle$ and $-1, \theta \in \mathcal{U}(\mathbb{Z}[\theta])\cap \langle \mathbb{S}_1\rangle$ so we have to prove $\mu_i \in \langle \mathbb{S}_1 \rangle$ where $i\in I_{p^n}.$
If $i\in I_{p^n}$ then  $i\in \mathcal{U}(\mathbb{Z}_{p^n})=\langle \bar{t} \rangle$ and we have $i \equiv t^j \ (\mbox{mod} \ p^n)$, for some integer $j$ such that $1 \leq j \leq \phi(p^n).$
Thus  $\mu_i=\mu_{t^j}.$ 
Let $\mu_i \in \mathcal{U}(\mathbb{Z}[\theta])$, where $i\in I_{p^n}$. We will prove that $\mu_i=\mu_{t^j}\in \langle \mathbb{S}_1\rangle$ in the following cases:
\begin{itemize}
\item[1)] $1\leq j \leq \frac{\phi(p^n)}{2}-1;$
\item[2)] $\frac{\phi(p^n)}{2} < j  < \phi(p^n);$
\item[3)] $j=\frac{\phi( p^n)}{2};$
\item[4)] $j=\phi(p^n)$.
\end{itemize}

The first case is immediate by definition of $\mathbb{S}_1$.
 If $\frac{\phi(p^n)}{2} < j  < \phi(p^n)$,  denote $r=j-\frac{\phi(p^n)}{2}.$ Then $t^j\equiv t^{r+\frac{\phi(p^n)}{2}} \ (\mbox{mod} \ p^n).$ 
As  $i\equiv t^j\equiv t^{r+\frac{\phi(p^n)}{2}} \ (\mbox{mod} \ p^n)$ and $t^{\frac{\phi(p^n)}{2}} \equiv -1 \ (\mbox{mod} \ p^n)$ it follows that  $i+t^r\equiv 0 \ (\mbox{mod} \ p^n)$.Thus $\mu_i=-\theta^i\mu_{t^r}$.
Since $r=j-\frac{\phi(p^n)}{2}$ and $\frac{\phi(p^n)}{2} < j  < \phi(p^n),$ we have  $1\leq r \leq \frac{\phi(p^n)}{2}-1$ and then $\mu_{t^r} \in \mathbb{S}_1$. Hence $\mu_i=-\theta^i\mu_{t^r}\in \langle \mathbb{S}_1\rangle.$
If $j=\frac{\phi(p^n)}{2},$ then  $i\equiv t^{\frac{\phi(p^n)}{2}} \ (\mbox{mod} \ p^n)$.
As $t^{\frac{\phi(p^n)}{2}} \equiv -1 \ (\mbox{mod} \ p^n)$ and $-1\equiv (p^n-1) \ (\mbox{mod} \ p^n),$ then  $i\equiv (p^n-1) \ (\mbox{mod} \ p^n)$.
 Therefore $\mu_i=\mu_{p^n-1}=1+\theta+\theta^2+\cdots+\theta^{p^n-2}=-\theta^{p^n-1} \in \langle \mathbb{S}_1\rangle$.
In the case $j=\phi(p^n)$ we have $i\equiv t^{\phi(p^n)} \ (\mbox{mod} \ p^n)$. Since $t^{\phi(p^n)}\equiv 1 \ (\mbox{mod} \ p^n)$ and  $-1\in \langle \mathbb{S}_1\rangle$ we conclude that $\mu_i=\mu_1=1 \in \langle \mathbb{S}_1\rangle.$ 
Consequently  $\mathcal{U}(\mathbb{Z}[\theta])\subseteq \langle \mathbb{S}_1 \rangle.$ $\square$

\

From now on $\kappa$ will denote the number $\frac{\phi(p^n)}{2}-1.$
As on the paper \cite{Raul} we change suitably the set of generators. The following proposition plays  the same role that Proposition  2.2 plays in \cite{Raul}.

\begin{prop1}
Let $t$ be an integer such that $\langle \bar{t} \rangle =\mathcal{U}(\mathbb{Z}_{p^n})$ and $$\mathbb{S}_2=\{-1, \theta, \mu_t, \mu_t^{-2}\mu_{t^2}, \mu_t^{-3}\mu_{t^3}, \cdots, \mu_t^{-\kappa}\mu_{t^{\kappa}}\}.$$ Then $\langle \mathbb{S}_1\rangle = \langle \mathbb{S}_2 \rangle.$
\end{prop1}

\prova

Clearly $\langle \mathbb{S}_2\rangle \subseteq \langle \mathbb{S}_1 \rangle,$ because $-1, \theta, \mu_{t^j} \in \langle \mathbb{S}_1\rangle,$ for all  $1 \leq j \leq \frac{\phi(p^n)}{2}-1$.
Let $\mu_{t^i} \in \mathbb{S}_1,$ then $\mu_{t^i}=\mu_t^{i}(\mu_t^{-i}\mu_{t^i}).$ Since $\mu_t^{i}$ and $(\mu_t^{-i}\mu_{t^i})$ are elements in $\langle \mathbb{S}_2\rangle,$ it follows that $\mu_{t^i}\in \langle \mathbb{S}_2\rangle.$ Hence $\langle \mathbb{S}_1\rangle = \langle \mathbb{S}_2 \rangle.$ $\square$

\

\begin{def1}
Let $\mathcal{U}(\mathbb{Z}_{p^n})=\langle \bar{t} \rangle.$ Denote by $\mathcal{U}$ the following subset of $\mathcal{U}(\mathbb{Z}[\theta])$:
$$\mathcal{U}:=  \{\theta, -(\mu_t^{\frac{p-1}{2}}), \mu_t^{-2}\mu_{t^2}, \cdots, \mu_t^{-\kappa}\mu_{t^{\kappa}} \}.$$
\end{def1}

Let $\Psi : \mathbb{Z}[\theta]\rightarrow \mathbb{Z}_{p}$ the  homomorphism defined by $\Psi(\theta)=\bar{1}$
and $\psi$ the restriction of $\Psi$ to $\mathcal{U}(\mathbb{Z}[\theta]).$ Clearly the image of $\psi$ is contained in $\mathcal{U}(\mathbb{Z}_{p}).$

\begin{teo1} \label{psi}
Let $\psi : \mathcal{U}(\mathbb{Z}[\theta]) \rightarrow \mathcal{U}(\mathbb{Z}_{p})$ the group homomorphism defined previously.
Then $ker(\psi)=\langle \mathcal{U}\rangle.$
\end{teo1}

\prova

Observe that $\psi(\mu_i)=i$ and if $t\in \mathbb{Z}$ such that  $t$ generates $\mathcal{U}(\mathbb{Z}_{p^n})$ than $t$ generates $\mathcal{U}(\mathbb{Z}_p).$ 
Applying $\psi$  for elements of $\mathcal{U},$ we have $\psi(\theta)=1,$
 $\psi(\mu_t^{-i}\mu_{t^i})=(t)^{-i}t^i=1,$ for all $2 \leq i \leq \kappa=\frac{\phi(p^n)}{2}-1$ and
 $\psi(-(\mu_t^{\frac{p-1}{2}}))=-(t^{\frac{p-1}{2}}) \equiv 1 \ (\mbox{mod}\ p),$ because $t^{\frac{p-1}{2}}\equiv -1 \ (\mbox{mod} \ p).$ Then $\langle \mathcal{U}\rangle \subseteq ker(\psi).$

Now, let $v\in ker(\psi) \subseteq \mathcal{U}(\mathbb{Z}[\theta])= \langle\mathbb{ S}_2\rangle$. Thus $v$ is of the form $$v=(-1)^{\alpha}\theta^{s_0}\mu_t^{s_1}(\mu_t^{-2}\mu_{t^2})^{s_2} \cdots (\mu_t^{-\kappa}\mu_{t^{\kappa}})^{s_{\kappa}},$$ where $0 \leq \alpha \leq 1,$ $ \ 0 \leq s_0 \leq p^n-1$ and $ \ s_j\in \mathbb{Z},$ for all $ j \in \{1,2, \cdots, \kappa\}$.
As we have $\theta,  \ (\mu_t^{-j}\mu_{t^j}) \in ker(\psi),$ then
$$v\in ker(\psi) \Leftrightarrow \left\lbrace
\begin{array}{lcl}
\alpha = 1 & \mbox{and} & t^{s_1}\equiv -1 \ (\mbox{mod} \ p)\\
\mbox{or} & \ & \\
\alpha =0 & \mbox{and} & t^{s_1}\equiv 1 \ (\mbox{mod} \ p)
\end{array}
\right.
$$
In the first case, since we have  $t^{\frac{p-1}{2}}\equiv -1 \ (\mbox{mod} \ p)$ then  $t^{s_1+\frac{(p-1)}{2}}\equiv 1 \ (\mbox{mod} \ p)$. As $\bar{t}$ generates $\mathcal{U}(\mathbb{Z}_p),$ we  shall have $s_1+\frac{(p-1)}{2}=(p-1)q,$ for some $q\in \mathbb{Z},$ and consequently $s_1=\frac{(p-1)}{2}(2q-1).$ Then $\mu_t^{s_1}=(\mu_t^{\frac{p-1}{2}})^{2q-1}.$ If $v\in ker(\psi),$ we have $$\begin{array}{lcl}
v&=&(-1)^{\alpha}\theta^{s_0}\mu_t^{s_1}(\mu_t^{-2}\mu_{t^2})^{s_2} \cdots (\mu_t^{-\kappa}\mu_{t^{\kappa}})^{s_{\kappa}}\\
&=&(1)\theta^{s_0}(-\mu_t^{\frac{p-1}{2}})^{2q-1}(\mu_t^{-2}\mu_{t^2})^{s_2} \cdots (\mu_t^{-\kappa}\mu_{t^{\kappa}})^{s_{\kappa}},
\end{array}$$ then $v\in \langle\mathcal{U}\rangle.$ In the other case, we can consider $s_1=(p-1)q=\frac{(p-1)}{2}2q,$ for some $q\in \mathbb{Z}$ and then $\mu_t^{s_1}=(\mu_{t}^{\frac{p-1}{2}})^{2q}=(-\mu_t^{\frac{p-1}{2}})^{2q}.$
Since $\alpha=0,$ we obtain
$$v=(-1)^{\alpha}\theta^{s_0}\mu_t^{s_1}(\mu_t^{-2}\mu_{t^2})^{s_2} \cdots (\mu_t^{-\kappa}\mu_{t^{\kappa}})^{s_{\kappa}}=\theta^{s_0}(-\mu_t^{\frac{p-1}{2}})^{2q}(\mu_t^{-2}\mu_{t^2})^{s_2} \cdots (\mu_t^{-\kappa}\mu_{t^{\kappa}})^{s_{\kappa}},$$ so $v \in \langle\mathcal{U}\rangle.$
Thus $ker(\psi)\subseteq \langle \mathcal{U}\rangle.$  $\square$

\

The next definition is similar to the definition on \cite{Raul}.

\begin{def1}
Given a positive integer $q$ and an integer $s$ relatively prime with $p^n$, we introduce an unit of $\mathbb{Z}[\theta]:$
$$\hspace{1cm} \omega_{q,s}=\displaystyle\sum_{j=0}^{q-1}\theta^{js}=1+\theta^s+\theta^{2s}+\cdots+\theta^{(q-1)s}. \hspace{2cm} (***)$$
\end{def1}
If $q$ is positive, then the unit $\mu_q$ is exactly the unit $\omega_{q,1}.$
This definition is valid for $p=2.$

\

The following proposition is similar to the Proposition 2.4 on \cite{Raul} and the proofs are quite similar. For this reason it will be omitted.

\begin{prop1}  \label{omega}
Let $q$ and $s$ positive integers with $q$ relatively prime with $p^n.$ \linebreak According to the notation above, the following equality holds
$$\mu_{q^s}=\displaystyle\prod_{j=0}^{s-1}\omega_{q,q^j}=\omega_{q,1}\omega_{q,q}\omega_{q,q^2}\cdots\omega_{q,q^{s-1}}$$
\end{prop1}

%
%
%

\begin{def1}
Let $t$ be an integer such that $\bar{t}$ generates  $\mathcal{U}(\mathbb{Z}_{p^n})$. Define $h_i$ as the following unit of $\mathcal{U}(\mathbb{Z}[\theta]):$ $$h_i=\omega_{t,1}^{-1}\omega_{t,t^i}.$$
\end{def1}

\begin{teo1} \label{ultimo}
Let $\mathcal{U}_0=\{\theta, h_1,  h_2,  \cdots, \ h_{\kappa}\}$ and $\mathcal{U}'=\{\theta, \ \mu_t^{-2}\mu_{t^2} \ ,\cdots, \ \mu_{t}^{-(\kappa+1)}\mu_{t^{\kappa+1}}\}.$  Then $\langle\mathcal{U}_0\rangle =\langle\mathcal{U}'\rangle.$
\end{teo1}

\prova

By the Proposition \ref{omega}, we have  $$\mu_{t^i}=\omega_{t,1}\omega_{t,t}\omega_{t,t^2}\cdots\omega_{t,t^{i-1}}.$$
Since $\omega_{t,1}=\mu_t,$ it follow that  $$\mu_{t^i}=\mu_t\omega_{t,t}\omega_{t,t^2}\cdots\omega_{t,t^{i-1}}$$ and multiplying the both sides of this equality by  $\mu_t^{-i},$ we yield
\begin{eqnarray*}
\mu_t^{-i}\mu_{t^i}&=&\mu_t^{-i+1}(\omega_{t,t}\omega_{t,t^2}\cdots\omega_{t,t^{i-1}})\\
&=&(\mu_t^{-1}\omega_{t,t})(\mu_t^{-1}\omega_{t,t^2})\cdots(\mu_t^{-1}\omega_{t,t^{i-1}})\\
&=&h_1h_2\cdots h_{i-1} \in \langle \mathcal{U}_0 \rangle ,
\end{eqnarray*}
for all integer $i,$  $2 \leq i \leq \kappa+1.$ Therefore $\mathcal{U}' \subseteq \langle \mathcal{U}_0 \rangle.$
On the other hand, by the Proposition above, we have $$\mu_{t^{i+1}}=\underbrace{\omega_{t,1}\omega_{t,t}\omega_{t,t^2}\cdots\omega_{t,t^{i-1}}}_{\mu_{t^i}}\omega_{t,t^i}=\mu_{t^i}\omega_{t,t^i}.$$
From this equality we obtain $\omega_{t,t^i}=\mu_{t^{i+1}}\mu_{t^i}^{-1}.$ Then
\begin{eqnarray*}
h_i&=&\omega_{t,1}^{-1}\omega_{t,t^i}=\omega_{t,1}^{-1}\mu_{t^{i+1}}\mu_{t^i}^{-1}=\mu_t^{-1}\mu_t^i\mu_t^{-i}\mu_{t^{i+1}}\mu_{t^i}^{-1}\\
&=&(\mu_t^{-(i+1)}\mu_{t^{i+1}})(\mu_t^{-i}\mu_{t^i})^{-1} \in \langle \mathcal{U}' \rangle .
\end{eqnarray*}
Therefore, $h_i \in \langle \mathcal{U}' \rangle$, for all integer $i$, $1\leq i \leq \kappa.$
It proves  $\mathcal{U}_0\subseteq \langle\mathcal{U}'\rangle.$ $\square$

\

\section{When $p=2$}

When $n\geq 3$ it is known that $\mathcal{U}(\mathbb{Z}_{2^n})$ is not cyclic, $\mathcal{U}(\mathbb{Z}_{2^n})=\langle \bar{3}\rangle \times \langle \overline{-1} \rangle$  and the order of $\bar{3}$ is $2^{n-2}.$ These facts suggest  change $t$ in the previous section by $3$ and it can be proved easily in the next proposition.  

\begin{prop1} \label{potenciade3}
Let $\theta$  be the $2^n$--th primitive root of unity and $i=( 2^{n}-1)^s3^q,$ where $0 \leq s \leq 1$ and $0 \leq q \leq 2^{n-2}-1.$ The following statements hold:
\begin{itemize}
\item[(i)]  $s=0  \Longrightarrow  \mu_i=\mu_{3^q};$
\item[(ii)] $s=1 \Longrightarrow \mu_i=-\theta^i\mu_{3^q}$
\end{itemize}
\end{prop1}

\prova

The proof of {\it(i)} is trivial.
If $s=1$ and $ q\neq 0,$ we have $ 1\leq i=(2^{n}-1)3^q \leq 2^n3^q $. Moreover $2^n3^q=i+3^q,$ and then $i+3^q\equiv 0 \ (\mbox{mod} \ 2^n).$ Hence  $\mu_i=-\theta^i\mu_{3^q}.$  $\square$

\

The next two propositions are the same as in the previous section with a slight change in the proof.

 \begin{teo1}
 Let $i\in I_{2^n}=\{ i\in \mathbb{Z} \mid \ \ 1< i< 2^{n-1},  \ gcd(2^n,i)=1 \}$, $\theta$ be the $2^n$--th primitive root of unity and  $n\geq 3$   such that $\mathbb{S}=\{ \ \theta, \  \mu_i= 1+\theta+\cdots+\theta^{i-1} \mid   i\in I_{2^n}\}$ generates the group of units of $\mathbb{Z}[\theta].$ Then $\mathbb{S} _1=\{  \theta, \ \mu_3, \ \mu_{3^2}, \cdots, \ \mu_{3^{\kappa}} \}$ generates $\mathcal{U}(\mathbb{Z}[\theta]),$ where $\kappa=\frac{\phi(2^n)}{2}-1.$
 \end{teo1}

\prova

Clearly $\mathbb{S}_1\subseteq \mathcal{U}(\mathbb{Z}[\theta]).$ 
In the other case, since $\mathcal{U}(\mathbb{Z}[\theta])=\langle \{-1,\ \theta, \ \mu_i \ | \ i\in \mathcal{I}_{2^n}\}\rangle$ we have to verify $\mu_i\in \langle \mathbb{S}_1\rangle$ for $i\in \mathcal{I}_{2^n}.$
If $i\in \mathcal{I}_{2^n},$ it follows that \linebreak ${i} \in \mathcal{U}(\mathbb{Z}_{2^n})=\langle 2^{n}-1 \rangle \times \langle 3 \rangle$.
So we obtain $i=( 2^{n}-1)^s3^q,$ where $0 \leq s \leq 1$ and $0 \leq q \leq 2^{n-2}-1.$  By the Proposition \ref{potenciade3} we get $\mu_i=\mu_{3^q} \in \mathbb{S}_1$ or $\mu_i=-\theta^i\mu_{3^q} \in \langle \mathbb{S}_1 \rangle.$
Hence $\mu_i \in \langle \mathbb{S}_1 \rangle,$ for all $i\in \mathcal{I}_{2^n}.$ $\square$

\

\begin{teo1}
In the same hypothesis as the previous theorem, let 
 $$\mathbb{S}_2=\{  \theta, \ \mu_3, \ \mu_3^{-2}\mu_{3^2}, \ \mu_3^{-3}\mu_{3^3}, \cdots, \ \mu_3^{-\kappa}\mu_{3^{\kappa}} \},$$ where $\kappa=\frac{\phi(2^n)}{2}-1.$ Then
$\langle \mathbb{S}_2 \rangle = \langle \mathbb{S}_1 \rangle.$
\end{teo1}

\prova

Since $\mu_{3^i} \in \mathbb{S}_1,$ for all $1 \leq i \leq \kappa,$ it follows that $\mu_3^{-i}\mu_{3^i}\in \langle \mathbb{S}_1 \rangle.$ Moreover, $\theta\in \langle \mathbb{S}_1 \rangle \cap \langle \mathbb{S}_2 \rangle,$ thus $\mathbb{S}_2 \subseteq \langle \mathbb{ S}_1 \rangle.$
On the other hand, as $\mu_{3^i}\in \mathbb{S}_1$  and $\mu_3^i, \mu_3^{-i}\mu_{3^i}\in \mathbb{S}_2$ it yields $\mu_{3^i}=\mu_3^i(\mu_3^{-i}\mu_{3^i}) \in \langle \mathbb{S}_2 \rangle,$ with $1 \leq i \leq \kappa.$ 
So  $\langle \mathbb{S}_2 \rangle = \langle \mathbb{S}_1 \rangle.$ $\square$

\

When $p=2,$ the homomorphism $\psi$ defined in the preceding section is not necessary and we continue with the following.

\begin{prop1} \label{indice2n}
Let $\mathcal{U}'=\{\theta, \ \mu_3^{-2}\mu_{3^2}, \ \mu_3^{-3}\mu_{3^3}, \cdots, \ \mu_3^{-\kappa}\mu_{3^\kappa}, \ \mu_3^{-(\kappa+1)}\}$ a subset of $\mathcal{U}(\mathbb{Z}[\theta]).$ Then $\langle \mathcal{U}' \rangle$ is a subgroup of $\mathcal{U}(\mathbb{Z}[\theta])$ of rank $\kappa+1=\frac{\phi(2^n)}{2}=2^{n-2}.$
\end{prop1}

\prova

Observe that $2^{n-2}=\frac{\phi(2^n)}{2}=\kappa +1$ and $3^{\kappa +1}\equiv 1 \ (\mbox{mod} \ 2^n).$
Thus $\mu_{3^{\kappa+1}}=\mu_1=1$ and then $\mu_{3}^{-(\kappa+1)}\mu_{3^{\kappa +1}}=\mu_{3}^{-(\kappa+1)}.$ 

Since $\mathcal{U}(\mathbb{Z}[\theta])=\langle \mathbb{S}_2 \rangle = \langle  \{  \theta, \ \mu_3, \ \mu_3^{-2}\mu_{3^2}, \ \mu_3^{-3}\mu_{3^3} \cdots, \ \mu_3^{-\kappa}\mu_{3^{\kappa}} \} \rangle,$ we conclude that $|\mathcal{U}_1(\mathbb{Z}[\theta]):\langle \mathcal{U}' \rangle|=\kappa+1=\frac{\phi(2^n)}{2}.$ $\square$

\

Recall  $$\omega_{q,s}=1+\theta^s+\theta^{2s}+\cdots+\theta^{(q-1)s},$$ and let us define the following unit of $\mathcal{U}(\mathbb{Z}[\theta]):$
$$h_i=\omega_{3,1}^{-1}\omega_{3,3^i}, \ \ 1\leq i \leq \kappa.$$

Similarly to the Proposition \ref{ultimo}, we prove the proposition below.

\begin{teo1} \label{geradorV}
With the notation above, the set $\mathcal{U}_0=\{\theta, \ h_1, \ h_2,  \cdots, \ h_{\kappa} \}$ generates $\langle\mathcal{U}'\rangle.$
\end{teo1}

\section{Units of $\mathbb{Z}C_{p^n}$}

Let $C_{p^n}=\langle g \rangle$ be the cyclic group of order $p^n$ and $\pi_1:\mathbb{Z}C_{p^n}\rightarrow \mathbb{Z}[\theta]$ be the homomorphism defined by $\pi_1\left(\sum a_gg\right)=\sum a_g\theta.$  Denote by $\overline{\pi}_1$ the restriction of $\pi_1$ to $\mathcal{U}_1(\mathbb{Z}C_{p^n})$ and clearly the image of $\overline{\pi}_1$ is contained in $\mathcal{U}_1(\mathbb{Z}[\theta]),$ where $\mathcal{U}_1(\mathbb{Z}[\theta])$ is  $ker(\psi)$ when $p$ is an odd prime and $\mathcal{U}(\mathbb{Z}[\theta])$ when $p=2.$
The following result give  a general  description to $\mathcal{U}_1(\mathbb{Z}C_{p^n}).$

\

\begin{teo1} \label{geradoresdeU1ZCpn}
$\mathcal{U}_1(\mathbb{Z}C_{p^n})\cong ker(\overline{\pi}_1) \times  Im(\overline{\pi}_1).$
\end{teo1}

\prova

Let $\iota: ker(\overline{\pi}_1)\rightarrow \mathcal{U}_1(\mathbb{Z}C_{p^n})$ be the inclusion and consider the following short exact sequence
$$\{1\} \longrightarrow  ker(\overline{\pi}_1) \stackrel{\iota}{ \longrightarrow} \mathcal{U}_1(\mathbb{Z}C_{p^n}) \stackrel{\overline{\pi}_1}{ \longrightarrow} Im(\overline{\pi}_1)\longrightarrow \{1\}$$
If this sequence splits, then  $\mathcal{U}_1(\mathbb{Z}C_{p^n})\simeq ker(\overline{\pi}_1)\times Im(\overline{\pi}_1),$ proving the result.  

By the  Dirichlet unit theorem we have $\mathcal{U}_1(\mathbb{Z}[\theta])=\langle \theta \rangle \times L,$ where $L$ is a free abelian group of finite rank. Then  $Im(\overline{\pi}_1)=T\times L_1,$ where $T$  is a torsion group 
  and $L_1$ is a free abelian group. Since $\langle \theta \rangle=\overline{\pi}_1(\langle g \rangle ),$ we have $T=\langle \theta \rangle$ and then $Im(\overline{\pi}_1)=\langle \theta \rangle\times L_1.$

The set $L_1$ is free, so it is projective. Then there exists a homomorphism  \linebreak $\tau_1: L_1 \rightarrow \mathcal{U}_1(\mathbb{Z}C_{p^n})$ such that $\overline{\pi}_1\circ \tau_1(u)=u,$ for all $ u\in L_1.$ We define \linebreak  $\tau:Im(\overline{\pi}_1) \rightarrow \mathcal{U}_1(\mathbb{Z}C_{p^n})$ by $\tau(u\cdot \theta^i)=\tau_1(u)\cdot g^i.$ Thus $\tau $ is well defined and $\overline{\pi}_1 \circ \tau $ is the identity of ${Im(\overline{\pi}_1)}.$ 
Then the sequence splits. $\square$

\

Now we shall describe a set of independent generators to $ker(\overline{\pi}_1)$ and to an isomorphic subgroup of 
$Im(\overline{\pi}_1)$ in $\mathcal{U}_1(\mathbb{Z}C_{p^n}).$

\begin{prop1}  \label{continencia}
Let $ \ \mathcal{U}_0=\{ \theta, \  h_1, \ h_2,  \cdots, \ h_{\kappa}\}$, where $h_i$ is an unit of $\mathbb{Z}[\theta]$ defined by $$h_i =  \omega_{t,1}^{-1}\omega_{t,t^i}=(1+\theta+\theta^2+\cdots+\theta^{t-1})^{-1}(1+\theta^{t^i}+\theta^{2t^i}+\cdots+\theta^{(t-1)t^i}).$$  
\begin{enumerate}
\item  If $p$  is an odd prime and $t$ is an integer such that $\bar{t}$ generates  $\mathcal{U}(\mathbb{Z}_p),$ or
\item if $p=2$ and $t=3,$ 
\end{enumerate}
then $\langle\mathcal{U}_0\rangle \subseteq Im(\overline{\pi}_1).$
\end{prop1}

\prova

Clearly $\theta\in Im(\overline{\pi}_1),$ so
we shall prove that $h_i \in Im(\overline{\pi}_1), \ 1\leq i \leq \kappa.$ 
\begin{eqnarray*}
h_i & = & \omega_{t,1}^{-1}\omega_{t,t^i}=(1+\theta+\theta^2+\cdots+\theta^{t-1})^{-1}(1+\theta^{t^i}+\theta^{2t^i}+\cdots+\theta^{(t-1)t^i}) \\
& = & (1+\theta^t+\theta^{2t}+\cdots+\theta^{(r-1)t})(1+\theta^{t^i}+\theta^{2t^i}+\cdots+\theta^{(t-1)t^i}),
\end{eqnarray*}
 where $r$ is the least integer such that $tr\equiv 1 \ (\mbox{mod} \ p^n).$

Observe that $${\vartheta}_i=(1+g^t+g^{2t}+\cdots+ g^{(r-1)t})(1+g^{t^i}+g^{2t^i}+\cdots+g^{(t-1)t^i}) - \frac{(tr-1)}{p^n}\widehat{g} $$ is the Hoechsmann unit in $\mathbb{Z}C_{p^n},$ where $\widehat{g}=1+g+g^2+\cdots +g^{p^n-1}.$
Since $\overline{\pi}_1({\vartheta}_i)=h_i,$ it follows $\mathcal{U}_0\subseteq Im(\overline{\pi}_1).$ $ \square$

\

With a suitable hypothesis the equality holds.

\begin{prop1} 
If $(-1)^p(\mu_t)^{\frac{\kappa+1}{p}}=(-1)^p(\mu_{t})^{\frac{\phi(p^n)}{2p} }\not\in Im(\overline{\pi}_1),$ then
$$|\mathcal{U}_1(\mathbb{Z}[\theta]):Im(\overline{\pi}_1)|= \left\{ 
\begin{array}{rcl} 
p^{n-1},& \mbox{if} &  p \mbox{ \ is odd prime}\\ 
p^{n-2}, & \mbox{if} & p=2
\end{array} 
\right. $$ Moreover, $Im(\overline{\pi}_1)=\langle \mathcal{U}_0 \rangle .$
\end{prop1}

\prova

We have the following chain
$$\langle\mathcal{U}_0\rangle \subseteq Im(\overline{\pi}_1)\subsetneq \mathcal{U}_1(\mathbb{Z}[\theta]).$$
If $p$ is odd prime, we have $\mathcal{U}_1(\mathbb{Z}[\theta])=\langle \{\theta, -(\mu_t^{\frac{p-1}{2}}), \mu_t^{-2}\mu_{t^2} \ ,\cdots, \ \mu_{t}^{-\kappa}\mu_{t^{\kappa}}\}\rangle,$
   and  $\langle \mathcal{U}' \rangle =\langle\{\theta, \ \mu_t^{-2}\mu_{t^2} \ ,\cdots, \ \mu_{t}^{-(\kappa+1)}\mu_{t^{\kappa+1}}\}\rangle =\langle \mathcal{U}_0 \rangle.$ 

Observe that 
$\mu_t^{-(\kappa+1)}\mu_{t^{\kappa+1}}=\theta^{-1}\cdot (-\mu_t^{-\frac{p-1}{2}})^{p^{n-1}}\in Im(\overline{\pi}_1)$ and $\theta^{-1}\in Im(\overline{\pi}_1),$ \linebreak   then $(-\mu_t^{-\frac{p-1}{2}})^{p^{n-1}}\in Im(\overline{\pi}_1).$ By the hypothesis,
 $(-1)^p(\mu_t)^{\frac{\kappa+1}{p}}\not\in Im(\overline{\pi}_1),$  
so  $|\mathcal{U}_1(\mathbb{Z}[\theta]):Im(\overline{\pi}_1)|=p^{n-1}.$ 
Since $|\mathcal{U}_1(\mathbb{Z}[\theta]):\mathcal{U}_0|=p^{n-1},$ we have   $Im(\overline{\pi}_1)=\langle \mathcal{U}_0 \rangle .$

When $p=2,$ we have already proved that  $$|\mathcal{U}_1(\mathbb{Z}[\theta]):\langle \mathcal{U}' \rangle|=|\mathcal{U}_1(\mathbb{Z}[\theta]):\langle \mathcal{U}_0 \rangle|=\kappa +1=\frac{\phi(2^n)}{2}=2^{n-2}=p^{n-2}.$$
 As $\mu_3^{2^{n-2}}=(\mu_{3^{\kappa+1}}\mu_3^{-(\kappa+1)})^{-1}\in Im(\overline{\pi}_1)$ and, by the hypothesis, $\mu_3^{2^{n-3}}\not\in Im(\overline{\pi}_1),$ thus  $|\mathcal{U}_1(\mathbb{Z}[\theta]):Im(\overline{\pi}_1)|=2^{n-2}.$ Hence  $Im(\overline{\pi}_1)=\langle \mathcal{U}_0 \rangle .$ $\square$
 
\

\begin{teo1} \label{iso}
Let $S=\{ \ g, \ {\vartheta}_1, \ {\vartheta}_2,  \ \vartheta_3, \  \cdots \ , \ {\vartheta}_{\kappa}\}$ a subset of $\mathcal{U}_1(\mathbb{Z}C_{p^n}),$ where \linebreak ${\vartheta}_i=(1+g^t+g^{2t}+\cdots+ g^{(r-1)t})(1+g^{t^i}+g^{2t^i}+\cdots+g^{(t-1)t^i}) - \frac{(tr-1)}{p^n}\hat{g} \in \mathcal{U}_1(\mathbb{Z}C_{p^n})$ is the Hoechsmann unit and consider  $(-1)^p(\mu_{t})^{\frac{\phi(p^n)}{2p} }\not\in Im(\overline{\pi}_1).$ Then
$\langle S \rangle \cong \langle \mathcal{U}_0 \rangle.$
\end{teo1}

\prova

The homomorphism $\overline{\pi}_1:\mathcal{U}_1(\mathbb{Z}C_{p^n}) \rightarrow Im(\overline{\pi}_1)=\langle \mathcal{U}_0 \rangle$ is surjective and \linebreak $\overline{\pi}_1(S)=\mathcal{U}_0.$
We will show that $\overline{\pi}_1|_{\langle S \rangle}: \langle S \rangle \rightarrow \langle \mathcal{U}_0 \rangle$ is injective.

Let $\overline{\pi}_1(g^{i_0}{\vartheta}_1^{i_1}  {\vartheta}_2^{i_2} \cdots \ {\vartheta}_{\kappa}^{i_{\kappa}})=1,$ where $ 0\leq i_0 \leq p^n-1$ and $i_1, \cdots, i_{\kappa}\in \mathbb{Z}.$ 
 Then $\theta^{i_0}h_1^{i_1}  h_2^{i_2} \cdots h_{\kappa}^{i_{\kappa}}=1,$ thus  $1=(\theta^{i_0}h_1^{i_1} h_2^{i_2} \cdots h_{\kappa}^{i_{\kappa}})^{p^n}=h_1^{i_1p^n}h_2^{i_2p^n}\cdots h_{\kappa}^{i_{\kappa}p^n}.$ The elements $h_1, h_2, \cdots, h_{\kappa}$ are multiplicatively independent, therefore $i_j=0,$  $1\leq j \leq {\kappa}.$ Then $1=\overline{\pi}_1(x)= \overline{\pi}_1(g^{i_0})=\theta^{i_0}$ and consequently $i_0=0.$ Thus  $\langle S \rangle \cong \langle \mathcal{U}_0 \rangle.$  $\square$
 
\

The hypothesis $(-1)^p(\mu_{t})^{\frac{\phi(p^n)}{2p} }\not\in Im(\overline{\pi}_1)$ is valid in each of the cases when $\phi(p^n)\leq 66.$ 
We restricted for these powers of primes because we could not find a general way to prove it and in these cases  it is known that the set $\mathbb{S}$ generates $\mathcal{U}(\mathbb{Z}[\theta]).$
The proof of this hypothesis requires lots of calculations as well as the way we characterized $ ker(\overline{\pi}_1).$ We will comment more about this hypothesis later.  
Now we will analyse the other component of $\mathcal{U}_1(\mathbb{Z}C_{p^n})$ we are looking for.

\

\begin{prop1} \label{simetriadonucleo}
 The elements in  $ ker(\overline{\pi}_1)$ are normalized symmetric units.
\end{prop1}

\prova

Let  $w\in ker(\overline{\pi}_1)$ and $w=g^i\cdot w',$ for some integer $i,$ $0\leq i\leq p^n-1$ and \linebreak $w'\in \mathcal{U}_{*}(\mathbb{Z}C_{p^{n}}).$ Then $1=\overline{\pi}_1(w)=\theta^i \cdot \overline{\pi}_1(w')$ and then $\overline{\pi}_1(w')=\theta^{-i}.$ The homomorphism preserves symmetry, so  $\overline{\pi}_1(w')=\theta^{-i}$ is symmetric because $w'$ is.
When $p$ is odd prime, the unique integer $i$ such that $\theta^{-i}$ is symmetric is $i=0.$ Then $w=g^0\cdot w'=w'$ is symmetric unit. In the case $p=2,$ the element  $\theta^{-i}$ is symmetric when $i=0$ or $i=2^{n-1}.$ In these cases we have $w=g^0 \cdot w'=w'$ or $w=g^{2^{n-1}}\cdot w'.$ Since $g^{2^{n-1}}$ and $w'$ are symmetric then $w$ is symmetric unit. $\square$

\

Before explain the method we used to find the generator of $ ker(\overline{\pi}_1)$, let us see a result in modular group rings and recall a general structure about units of integral group rings.

\begin{def1}
Let $\mathbb{F}G$ denote the group ring of a finite $p$-group $G$ over the prime field $\mathbb{F}.$ We define the group of normalized units in $\mathbb{F}G$ by  $$\mathcal{U}_1(\mathbb{F}G)=\left\lbrace\displaystyle\sum_{x\in G}a_xx\in \mathcal{U}(\mathbb{F}G) \ | \ \sum a_x\equiv 1 \ (\mbox{mod} \ p) \right\rbrace .$$
\end{def1}

It is known that $|\mathcal{U}_1(\mathbb{F}G)|=p^{|G|-1}$ and that $\mathcal{U}_1(\mathbb{F}G)$ is a normal Sylow \linebreak $p$-subgroup of the full group of units in $\mathbb{F}G.$

\begin{def1}
A set $\{u_1, u_2, \cdots, u_k\}$ is multiplicatively independent in a finite abelian $p$-group  if
$$\langle u_i \rangle \cap \langle u_1, u_2, \cdots , u_{i-1}, u_{i+1}, \cdots, u_k \rangle =\{1\}, \ \ \mbox{for all} \ 1\leq i \leq k,$$
i.e., if $u_1^{r_1}u_1^{r_2}\cdots u_k^{r_k}=1$ then $u_j^{r_j}=1,$ for all $1\leq j \leq k.$
\end{def1}

\begin{prop1} \label{independente0}
Let $C_{p^n}=\langle g \rangle, \ x=g-1$ and $u,v\in \mathcal{U}_1(\mathbb{Z}_pC_{p^n})$ given by:
$$u=1+a_0x^{p^iq}+a_1x^{p^iq+1}+a_2x^{p^iq+2}+\cdots +a_tx^{p^n-1}, \ \ t=p^n-1-p^iq,$$
$$v=1+b_0x^{p^jr}+b_1x^{p^jr+1}+b_2x^{p^jr+2}+\cdots +b_{s}x^{p^n-1}, \  \ s=p^n-1-p^jr,$$
where $ a_0\not\equiv 0 (\mbox{mod} \ p),$  $b_0\not\equiv 0 (\mbox{mod} \ p),$  $q$ and $r$ distinct integers such that $p\nmid q$ and $p\nmid r.$  Then the set $\{u, v\}$ is multiplicatively independent.
\end{prop1}
\prova

Let us consider $u\neq 1$ and $v\neq 1.$ 
If $w\in \langle u \rangle \cap \langle v \rangle,$
then there exist integers  $n_1, n_2$ such that
$$w=u^{n_1}=(1+a_0x^{p^iq}+a_1x^{p^iq+1}+a_2x^{p^iq+2}+\cdots +a_tx^{p^n-1})^{n_1}$$
and
$$w=v^{n_2}=(1+b_0x^{p^jr}+b_1x^{p^jr+1}+b_2x^{p^jr+2}+\cdots +b_{s}x^{p^n-1})^{n_2}.$$
  Let $n_1=p^km_1$ and $n_2=p^lm_2,$ with $m_1$ and $m_2$ being integers such that  $p$ not divide neither $m_1$ nor $m_2.$ Therefore
\begin{eqnarray}
\nonumber w=u^{n_1}&=&(1+a_0x^{p^iq}+a_1x^{p^iq+1}+a_2x^{p^iq+2}+\cdots +a_tx^{p^n-1})^{n_1}\\
 \nonumber & =&(1+a_0^{p^k}x^{qp^{i+k}}+a_1^{p^k}x^{(qp^{i}+1)p^k}+a_2^{p^k}x^{(qp^{i}+2)p^k}+\cdots +a_t^{p^k}x^{(p^n-1
)p^k})^{m_1}\\
 \label{eq7} &=&1+m_1x^{qp^{i+k}}+{c}_1x^{qp^{i+k}+1}+{c}_2x^{qp^{i+k}+2}+\cdots +{c}_{t'}x^{p^n-1}.
\end{eqnarray}
On the other hand,  we have 
\begin{eqnarray}
\nonumber w=v^{n_2}&=&(1+b_0x^{p^jr}+b_1x^{p^jr+1}+b_2x^{p^jr+2}+\cdots +b_{s}x^{p^n-1})^{n_2}\\
\nonumber &=&(1+b_0^{p^l}x^{rp^{j+l}}+b_1^{p^l}x^{(p^jr+1)p^l}+b_2^{p^l}x^{(p^jr+2)p^l}+\cdots +b_{s}^{p^l}x^{(p^n-1)p^l})^{m_2}\\
 \label{eq8}&=&1+m_2x^{rp^{j+l}}+{\beta}_1x^{rp^{j+l}+1}+{\beta}_2x^{rp^{j+l}+2}+\cdots +{\beta}_{s'}x^{p^n-1}.
\end{eqnarray}
As $q\neq r,$  $p\nmid q$ and $p\nmid r$ thus   ${qp^{i+k}}\neq {rp^{j+l}}.$ Let us suppose by contradiction that  $qp^{i+k} < p^n$ and $rp^{j+l}<p^n.$
In this case, we have $x^{qp^{i+k}}\neq x^{rp^{j+l}}.$ Then 
$m_1x^{qp^{i+k}}=0$ or $m_2x^{rp^{j+l}}=0.$ As much as $p$ not divide neither $m_1$ nor $m_2,$   we shall have  $x^{qp^{i+k}}=0$ or $x^{rp^{j+l}}=0,$ and it is a contradiction.  Thus  $qp^{i+k} \geq  n$ or $rp^{j+l}\geq n.$  These cases implies  $x^{qp^{i+k}}=0$ or $x^{rp^{j+l}}=0.$ Hence $w=1$ and  $\langle u \rangle \cap \langle v \rangle=\{1\}.$ $\square$

\
 
Consider $C_{p^n}=\langle g \rangle$ and $C_{p^{n-1}}=\langle h \rangle.$ Define $\pi_2:\mathbb{Z}C_{p^n}\rightarrow \mathbb{Z}C_{p^{n-1}}$ by \linebreak $\pi_2\left( \sum a_ig^i\right)=\sum a_ih^i,$ $f_1:\mathbb{Z}C_{p^{n-1}}\rightarrow \mathbb{Z}_pC_{p^{n-1}}$ the natural homomorphism and $f_2:\mathbb{Z}[\theta]\rightarrow \mathbb{Z}_pC_{p^{n-1}}$ by $f_2\left( \sum a_i\theta^i\right)=\sum \bar{a}_ih^i.$ We will restrict each of the homomorphism above to subgroups in the following way:
$$\overline{\pi}_2:=\pi_2|_{\mathcal{U}_1(\mathbb{Z}C_{p^n})}, \ \bar{f}_1:=f_1|_{\mathcal{U}_1(\mathbb{Z}C_{p^{n-1}})} \   \mbox{and} \ \bar{f}_2:=f_2|_{\mathcal{U}_1(\mathbb{Z}[\theta])}.$$ Therefore we have the commutative diagram:
$$
\xymatrix{ 
\mathcal{U}_1(\mathbb{Z}C_{p^n} )\ar[r]^{\overline{\pi}_1} \ar[d]_{\overline{\pi}_2} & \mathcal{U}_1(\mathbb{Z}[\theta]) \ar[d]^{\bar{f}_2} \\ 
\mathcal{U}_1(\mathbb{Z}C_{p^{n-1}}) \ar[r]_{\bar{f}_1} &  \mathcal{U}_1(\mathbb{Z}_pC_{p^{n-1}} )\\ 
} 
$$
It is known that $ker(\pi_1)$ is an ideal of $\mathbb{Z}C_{p^n}$ generated by  $\displaystyle\sum_{i=0}^{p-1}g^{ip^{n-1}}$ (see \cite{Karpilovsky}). We will denote this sum by $(\widehat{g^{p^{n-1}}}).$ We have $$ker(\pi_1)=\{a_0(\widehat{g^{p^{n-1}}})+a_1g(\widehat{g^{p^{n-1}}})+\cdots+a_{p^{n-1}-1}g^{p^{n-1}-1}(\widehat{g^{p^{n-1}}}) \ | \ a_i \in \mathbb{Z}\}$$
and then
$$ker(\overline{\pi}_1)=\{1+a_0(\widehat{g^{p^{n-1}}})+a_1g(\widehat{g^{p^{n-1}}})+\cdots+a_{p^{n-1}-1}g^{p^{n-1}-1}(\widehat{g^{p^{n-1}}}) \ | \ a_i \in \mathbb{Z}\}.$$

It is easy to prove that the rank of    $ker(\overline{\pi}_1)$ and $ker(\bar{f}_1)$ are equal, and that $\overline{\pi}_2|_{ker(\overline{\pi}_1)}$ is injective. If we find the set of independent generator of $ker(\bar{f}_1)$ then we can find to  $ker(\overline{\pi}_1),$ once we know that an element of $ker(\overline{\pi}_1)$ is of the form $1+a_0(\widehat{g^{p^{n-1}}})+a_1g(\widehat{g^{p^{n-1}}})+\cdots+a_{p^{n-1}-1}g^{p^{n-1}-1}(\widehat{g^{p^{n-1}}})$ and $\overline{\pi}_2|_{ker(\overline{\pi}_1)}$ is injective.

Let us recall a general structure about $\mathcal{U}(\mathbb{Z}G).$

It is   known that $\mathcal{U}(\mathbb{Z}G)=\pm \mathcal{U}_1(\mathbb{Z}G).$ When $G$ is a finite abelian group, then  $\mathcal{U}_1(\mathbb{Z}G)=G\times \mathcal{U}_2(\mathbb{Z}G),$ where $\mathcal{U}_2(\mathbb{Z}G)=\{u\in \mathcal{U}(\mathbb{Z}G) : u\equiv 1  \ \mbox{mod} \ (\Delta G)^2\}$ and $\Delta G$ is the augmentation ideal of $\mathbb{Z}G.$ 

Let $*$ be the standard involution in $\mathbb{Z}G$ given by $ u=\displaystyle\sum_{g\in G} a_gg \mapsto u^*=\displaystyle\sum_{g\in G} a_gg^{-1}.$ If $u^*=u,$ we will say that $u$ is $*-$symmetric  and $\mathcal{U}_{*}(\mathbb{Z}G)=\{u\in \mathcal{U}(\mathbb{Z}G): \ u^*=u\}$ will denote the set of symmetric units. It is known that $\mathcal{U}_2(\mathbb{Z}G)\subseteq \mathcal{U}_{*}(\mathbb{Z}G),$ and when $|G|$ is odd  $\mathcal{U}_1(\mathbb{Z}G)=G\times \mathcal{U}_*(\mathbb{Z}G).$ 

We will use symmetric units and we will explain the method used in each of the cases when $\phi(p^n)\leq 66$ in the following example. All the other cases can be found in:

\begin{center}
http://paginapessoal.utfpr.edu.br/kitani/units/UnitsZCpn.pdf/view
  \end{center}

\begin{exe1}
The rank of $Im(\overline{\pi}_1)$ is $2$ as well as the rank of  $\mathcal{U}_{*}(\mathbb{Z}C_{3^2}),$ so  \linebreak $\mathcal{U}_{*}(\mathbb{Z}C_{3^2})=Im(\overline{\pi}_1)$ and  $\ker(\overline{\pi}_1)=\{1\}.$ 

$${\vartheta}_i=(1+g^t+g^{2t}+\cdots+ g^{(r-1)t})(1+g^{t^i}+g^{2t^i}+\cdots+g^{(t-1)t^i}) - \frac{(tr-1)}{p^n}\widehat{g},$$ where $\widehat{g}=(1+g+g^2+\cdots+g^{p^n-1}).$ Therefore  
\begin{eqnarray*}
\vartheta_1 &=&(1+g^2+g^4+g^6+g^8)(1+g^{2})-\widehat{g}\\
&=& (g^2-g^3+g^4-g^5+g^6-g^7+g^8), \\
&=&g^5(-1+g-g^2+g^3+g^6-g^7+g^8)\\
\vartheta_2 &=& (1+g^2+g^4+ g^{6}+g^8)(1+g^4) - \widehat{g}\\
&=&(g^4-g^5+g^6-g^7+g^8)\\
  &=& g^6(1-g+g^2+g^7-g^8)
\end{eqnarray*}
$\vartheta_1^{'}=(-1+g-g^2+g^3+g^6-g^7+g^8)$ and $  \vartheta_2^{'}=(1-g+g^2+g^7-g^8)$ are the symmetric units of $\mathbb{Z}C_{3^2}.$  Then $$\mathcal{U}(\mathbb{Z}C_{3^2})=\pm \mathcal{U}_1(\mathbb{Z}C_{3^2})=\pm C_{3^2}\times \mathcal{U}_{*}(\mathbb{Z}C_{3^2})=\pm C_{3^2}\times \langle \{\vartheta_1^{'}, \ \vartheta_2^{'}\}\rangle.$$
\end{exe1}

\

\begin{exe1} Let us find $ker(\overline{\pi}_1)$ in the case $\mathbb{Z}C_{3^3}.$

Let  $C_{3^3}=\langle g \rangle,$  $C_9=\langle h \rangle$ and $x=h-1.$
Analysing the following diagram, we will find a subgroup $N$ of $ker(\bar{f}_1)$ such that $|\mathcal{U}_{*}(\mathbb{Z}C_9):N|=|Im(\bar{f}_1)|.$

$$
\xymatrix{ 
\mathcal{U}_1(\mathbb{Z}C_{3^3}) \ar[r]^{\overline{\pi}_1} \ar[d]_{\overline{\pi}_2} & \mathcal{U}_1(\mathbb{Z}[\theta] )\ar[d]^{\bar{f}_2} \\ 
\mathcal{U}_1(\mathbb{Z}C_{3^2} )\ar[r]_{\bar{f}_1} &  \mathcal{U}_1(\mathbb{Z}_3C_{3^2}) \\ 
} 
$$

It is known that $$\mathcal{U}_{*}(\mathbb{Z}C_9)= \langle  \{ u_1= -1+h-h^2+h^3+h^6-h^7+h^8, \  u_2=1-h+h^2+h^7-h^8\} \rangle.$$ 
As $x=h-1$ then $h=x+1$ and:
\begin{eqnarray*}
\bar{f}_1(u_1)&=&-\bar{1}+h-\bar{1}h^2+h^3+h^6-\bar{1}h^7+h^8=\bar{1}+\bar{2} x^4+\bar{2} x^5+x^6+x^7+x^8,\\
\bar{f}_1(u_2)&=&\bar{1}-\bar{1}h+h^2+h^7-\bar{1}h^8=\bar{1}+x^4+x^5+\bar{2}  x^7+\bar{2} x^8.
\end{eqnarray*}
The set $\{\bar{f}_1(u_1), \bar{f}_1(u_2)\}$ is not necessary multiplicatively independent. We find a multiplicatively independent set that generates $Im(\bar{f}_1)$ in the following way:
\begin{eqnarray*}
\bar{f}_1(u_1)&=&\bar{1}+\bar{2} x^4+\bar{2} x^5+x^6+x^7+x^8\\
\bar{f}_1(u_1)\cdot \bar{f}_1(u_2)&=&\bar{1}+x^6+\bar{2} x^8.
\end{eqnarray*}
By the Theorem \ref{independente0}, $\bar{f}_1(u_1)$ and $\bar{f}_1(u_1)\cdot \bar{f}_1(u_2)$ are independents, then \linebreak  $Im(\bar{f}_1)=\langle \bar{f}_1(u_1)\rangle \times \langle \bar{f}_1(u_1)\cdot \bar{f}_1(u_2) \rangle.$ The order of $\bar{f}_1(u_1)$  and $\bar{f}_1(u_1)\cdot \bar{f}_1(u_2)$ are 3, therefore $|Im(\bar{f}_1)|=9.$ 

We have 
$$
\bar{f}_1(u_1^3)=(\bar{1}+\bar{2} x^4+\bar{2} x^5+x^6+x^7+x^8)^3=\bar{1}+(\bar{2})^3 x^{4\cdot 3}+(\bar{2})^3 x^{5\cdot 3}+x^{6\cdot 3}+x^{7\cdot 3}+x^{8\cdot 3}=\bar{1} $$
and
$$ \bar{f}_1(u_2^3)=(\bar{1}+x^4+x^5+\bar{2}  x^7+\bar{2} x^8)^3=\bar{1}+x^{4\cdot 3}+x^{5\cdot 3}+(\bar{2})^3  x^{7\cdot 3}+(\bar{2})^3 x^{8\cdot 3}=\bar{1},$$
once $x^j=\bar{0}$ for $j\geq 9.$
 Let $N=\langle u_1^3, u_2^3 \rangle,$ then  $N\subseteq ker(\bar{f}_1).$ As $\mathcal{U}_{*}(\mathbb{Z}C_9)=\langle u_1, u_2 \rangle,$ we have $|\mathcal{U}_{*}(\mathbb{Z}C_9):N|=3^2=|Im(\bar{f}_1)|$ and hence  $N=ker(\bar{f}_1).$

Now, if  we take an element in $ ker(\overline{\pi}_1),$ it is of the form $$w_i=1+a_0\widehat{(g^9)}+a_1g\widehat{(g^9)}+a_2g^2\widehat{(g^9)}+a_3g^3\widehat{(g^9)}
+a_4g^4\widehat{(g^9)}+\cdots+a_7g^7\widehat{(g^9)}+a_8g^8\widehat{(g^9)},$$ then $\pi_2(w_i)=1+3\cdot (a_0+a_1h+a_2h^2+\cdots +a_8h^8).$
On the other hand $\pi_2(w_i)=u_i^3$ and then we can find the independent generators of $ker(\overline{\pi}_1):$ 

\begin{eqnarray*}
u_1^3 &=& 1+3(-12+11h-9h^2+6h^3-2h^4-2h^5+6h^6-9h^7+11h^8)=\pi_2(w_1)\\
&=& 1+3(a_0+a_1h+a_2h^2+a_3h^3+\cdots +a_7h^7+a_8h^8).
\end{eqnarray*}
Thus $a_0=-12, \ a_1=a_8=11, \  a_2=a_7=-9, \ a_3=a_6=6$ and $a_4=a_5=-2,$ so
\begin{small}
$$w_1=1-12\widehat{(g^9)}+11g\widehat{(g^9)}-9g^2\widehat{(g^9)}+6g^3\widehat{(g^9)}
-2g^4\widehat{(g^9)}-2g^5\widehat{(g^9)}+6g^6\widehat{(g^9)}-9g^7\widehat{(g^9)}+11g^8\widehat{(g^9)}.$$
\end{small} 
\begin{eqnarray*}
u_2^3 &=& 1+3(6-6h+5h^2-3h^3+h^4+h^5-3h^6+5h^7-6h^8)=\pi_2(w_2)\\
&=& 1+3(a_0+a_1h+a_2h^2+a_3h^3+\cdots +a_7h^7+a_8h^8)
\end{eqnarray*}
 Then $a_0=6, \ a_1=a_8=-6, \ a_2=a_7=5, \  a_3=a_6=-3$ and $a_4=a_5=1,$ and 
 \begin{small}
$$w_2=1+6\widehat{(g^9)}-6g\widehat{(g^9)}+5g^2\widehat{(g^9)}-3g^3\widehat{(g^9)}
+1g^4\widehat{(g^9)}+1g^5\widehat{(g^9)}-3g^6\widehat{(g^9)}+5g^7\widehat{(g^9)}-6g^8\widehat{(g^9)} .$$ 
\end{small}
Then
$ker(\overline{\pi}_1)=\langle \{w_1, \ w_2\} \rangle $
\end{exe1}

We used symmetric units because it is convenient. 
When $p=2,$ instead of $\mathcal{U}_{*}(\mathbb{Z}G)$ we prefer  make use of $\mathcal{U}_2(\mathbb{Z}G).$

%
%
%

\section{Validity of Hypothesis of Theorem \ref{iso}}

The hypothesis $(-1)^p(\mu_{t})^{\frac{\phi(p^n)}{2p} }\not\in Im(\overline{\pi}_1)$ is valid when $\phi(p^n)\leq 66.$
We could not prove it in a general way. We have proved using the calculations used to find the set of independents generators of $ker(\overline{\pi}_1).$ The following proposition give us a normalized unit of $\mathbb{Z}C_{p^r}.$
 
\begin{prop1}
Let $p$ be prime  number and $C_{p^{n-1}}=\langle h \rangle$. Consider $\mathcal{U}(\mathbb{Z}_{p^n})=\langle \bar{t} \rangle$ if $p$ is odd and $t=3$ if $p=2.$ Then
$$\varpi=(-1)^p(1+h+h^2+\cdots + h^{t-1})^{\frac{\phi(p^{n-1})}{2}}-(-1)^p\lambda\cdot \widehat{h}\in \mathcal{U}_1(\mathbb{Z}C_{p^{n-1}})$$
where $\widehat{h}=(1+h+h^2+\cdots +h^{p^{n-1}-1})$ and $\lambda=\dfrac{t^{\frac{\phi(p^{n-1})}{2}}-(-1)^p}{p^{n-1}}.$
\end{prop1}

\prova

Notice that $t^{\frac{\phi(p^{n-1})}{2}}\equiv (-1)^p \ (\mbox{mod} \ p^{n-1})$ thus $\lambda=\dfrac{t^{\frac{\phi(p^{n-1})}{2}}-(-1)^p}{p^{n-1}}\in \mathbb{Z},$ so $\varpi \in \mathbb{Z}C_{p^{n-1}}.$

\vspace{0.3cm}

 Let $\omega=(-1)^p(1+h^t+h^{2t}+\cdots +h^{t(s-1)})^{\frac{\phi(p^{n-1})}{2}}-(-1)^p\lambda' \cdot \widehat{h}$
 where $s$ is the least integer such that $ts\equiv 1 \ (\mbox{mod} \ p^{n-1})$ and $\lambda'=\dfrac{s^{\frac{\phi(p^{n-1})}{2}}-(-1)^p}{p^{n-1}}.$ We shall prove $\omega=\varpi^{-1}.$  

 As $ts\equiv 1 \ (\mbox{mod} \ p^{n-1})$ it follows  $t^{\frac{\phi(p^{n-1})}{2}}\cdot s^{\frac{\phi(p^{n-1})}{2}}\equiv 1 \ (\mbox{mod} \ p^{n-1}).$ Also \linebreak
 $t^{\frac{\phi(p^{n-1})}{2}}\equiv (-1)^p \ (\mbox{mod} \ p^{n-1}),$ then $ s^{\frac{\phi(p^{n-1})}{2}}\equiv (-1)^p \ (\mbox{mod} \ p^{n-1})$ and $\omega \in \mathbb{Z}C_{p^{n-1}}.$
 
 Observe that

 \noindent $\varpi \cdot \omega =\\
 =\left[(-1)^p(1+h+\cdots + h^{t-1})^{\frac{\phi(p^{n-1})}{2}}\right]\left[ (-1)^p(1+h^t+\cdots +h^{t(s-1)})^{\frac{\phi(p^{n-1})}{2}} \right]+R\widehat{h}$

 \noindent for some $R\in \mathbb{Z}$
 
 Thus $\varpi \cdot \omega =\left( 1+\dfrac{st-1}{p^{n-1}}\widehat{h}\right)^{\frac{\phi(p^{n-1})}{2}}+R\widehat{h}=1+R_1\widehat{h}, \mbox{for some }  \ R_1\in \mathbb{Z}.$

Since $\varpi$ and $\omega$ have augmentation $1,$ then $\varpi \cdot \omega=1+R_1\widehat{h}$  have also augmentation $1.$ Therefore $R_1=0$ and we have proved $\varpi \cdot \omega=1.$ $\square$
 
 \
 
%

\begin{prop1} \label{feio}
Let $n \in \mathbb{N}, \ n\geq 2$ and $p$ be prime number such that $\phi(p^n)\leq 66,$ and $\theta$ be a $p^n$-th primitive root of the unity. Consider $t\in \mathbb{Z}$ such that $\mathcal{U}(\mathbb{Z}_{p^n})=\langle \bar{t} \rangle$ if $p$ is odd prime and $t=3$ if $p=2$, and $\mu_t=1+\theta+\theta^2+\cdots+\theta^{t-1}$ a cyclotomic unit of $\mathbb{Z}[\theta].$
Then $(-1)^p(\mu_{t})^{\frac{\phi(p^n)}{2p} }\not\in Im(\overline{\pi}_1).$ 
\end{prop1}

\prova

Suppose by contradiction that  $v=(-1)^p(\mu_{t})^{\frac{\phi(p^n)}{2p} }=(-1)^p(\mu_{t})^{\frac{\phi(p^{n-1})}{2} }\in Im(\overline{\pi}_1).$ Then there exists  $u\in \mathcal{U}_1(\mathbb{Z}C_{p^n})$ such that $\overline{\pi}_1(u)=(-1)^p(\mu_{t})^{\frac{\phi(p^{n-1})}{2}}=v.$ The following diagram is commutative and then
 $\bar{f}_2(v)=\bar{f}_2(\overline{\pi}_1(u))=\bar{f}_1(\overline{\pi}_2(u))$ is an element in $ Im(\bar{f}_2)\cap Im(\bar{f}_1).$ 

$$
\xymatrix{ 
\mathcal{U}_1(\mathbb{Z}C_{p^n} )\ar[r]^{\overline{\pi}_1} \ar[d]_{\overline{\pi}_2} & \mathcal{U}_1(\mathbb{Z}[\theta]) \ar[d]^{\bar{f}_2} \\ 
\mathcal{U}_1(\mathbb{Z}C_{p^{n-1}}) \ar[r]_{\bar{f}_1} &  \mathcal{U}_1(\mathbb{Z}_pC_{p^{n-1}} )\\ 
} 
$$
Since $\varpi=(-1)^p(1+h+h^2+\cdots + h^{t-1})^{\frac{\phi(p^{n-1})}{2}}-(-1)^p\lambda\cdot \widehat{h} \in \mathcal{U}_1(\mathbb{Z}C_{p^{n-1}}),$ where 
$$\lambda=\dfrac{t^{\frac{\phi(p^{n-1})}{2}}-(-1)^p}{p^{n-1}},$$ 
thus  
  $$\bar{f}_1(\varpi)={(\overline{-1})^p}(\bar{1}+h+h^2+\cdots + h^{t-1})^{\frac{\phi(p^{n-1})}{2}}-(\overline{-1})^p\bar{\lambda}\cdot \widehat{h}=\bar{f}_2(v)-(\overline{-1})^p\bar{\lambda} \cdot \widehat{h}.$$ 

We have also $\bar{f}_2(v)\cdot (\overline{\lambda} \cdot \widehat{h})=(\overline{-1})^p\bar{t}^{\frac{\phi(p^{n-1})}{2}}\cdot  \overline{\lambda} \cdot \widehat{h}= \overline{\lambda} \cdot \widehat{h},$ because \linebreak  $\bar{f}_2(v)\cdot \widehat{h}=(\overline{-1})^p\bar{t}^{\frac{\phi(p^{n-1})}{2}}\cdot \widehat{h}$  and 
$(-1)^pt^{\frac{\phi(p^{n-1})}{2}}\equiv 1 \ (\mbox{mod} \ p^{n-1}).$ Then

 $$\bar{f}_1(\varpi)=\bar{f}_2(v)-(\overline{-1})^p\bar{\lambda} \cdot \widehat{h}=\bar{f}_2(v)\cdot \left[\bar{1}-(\overline{-1})^p\bar{\lambda} \cdot \widehat{h}\right].$$
Since $\bar{f}_2(v)=\bar{f}_1({\pi}_2(u))\in Im(\bar{f}_1),$  we have $ \left(\bar{1}-(\overline{-1})^p\bar{\lambda} \cdot \widehat{h}\right)\in Im(\bar{f}_1).$

If $x=h-1$ then 

\vspace{0.2cm}
\noindent $ \bar{1}-(\overline{-1})^p\bar{\lambda} \cdot \widehat{h}= \bar{1}-(\overline{-1})^p\bar{\lambda} \cdot (1+h+\cdots +h^{p^{n-1}-1})= \bar{1}-(\overline{-1})^p\bar{\lambda} \cdot x^{p^{n-1}-1}\in Im(\bar{f}_1).$

We have  $Im(\bar{f}_1)=\langle e_1 \rangle \times \cdots \times\langle e_l \rangle$  and easily we verify using the evaluations to find the set of independent generators of $ker(\overline{\pi}_1)$ that $\bar{1}-(\overline{-1})^p\bar{\lambda} \cdot x^{p^{n-1}-1}$ is multiplicatively  independent with $\{e_1, e_2, \cdots, e_l\}$ ( see the evaluations in \linebreak
 http://paginapessoal.utfpr.edu.br/kitani/units/UnitsZCpn.pdf/view ).
\noindent It is a contradiction, therefore we can not assume $v=(-1)^p(\mu_{t})^{\frac{\phi(p^n)}{2p} }\in Im(\overline{\pi}_1).$ $\square$
  
\begin{exe1}  In the case $ZC_{3^3}$ 
 we have found in the previous example \linebreak  $Im(\bar{f}_1)=\langle \{ \bar{1}+\bar{2} x^4+\bar{2} x^5+x^6+x^7+x^8, \  \bar{1}+x^6+\bar{2} x^8 \}\rangle.$
 
 Using the Theorem \ref{independente0}, it is easy to see that $\bar{1}-(\overline{-1})^p\bar{\lambda} \cdot x^{p^{n-1}-1}=\bar{1}+\bar{1}x^8$ is multiplicatively independent with $ \bar{1}+\bar{2} x^4+\bar{2} x^5+x^6+x^7+x^8$ and $ \bar{1}+x^6+\bar{2} x^8.$ Then  $v=(-1)(\mu_{t})^{\frac{\phi(3^3)}{2\cdot 3} }\not\in Im(\overline{\pi}_1).$
   \end{exe1}  

All the other cases can be found in:
\begin{center}
http://paginapessoal.utfpr.edu.br/kitani/units/UnitsZCpn.pdf/view
  \end{center}
  
%
%
  

\begin{thebibliography}{9}                                                                                                

\bibitem{Aleev 2} R. Z. Allev and Z. Panina, {\it The Units of Cyclic Groups of Orders 7 and 9}, Russian Math. 43, (2000), 80--83.

\bibitem {AyoubAyoub} { R. G. Ayoub  and C. Ayoub},  \textit{On the group ring of a finite abelian group}, Bull. Austr. Math.
Soc. 1 (1969), 245-261. 


%
%
%
%
%
%

\bibitem{Raul} R. A. Ferraz,  {\it Units of $\mathbb{Z}C_{p}$}, Contemporary Mathematics - American Mathematical Society (2009),  v. 499, 107--119.

%
%
%
%

\bibitem{Higman} G. Higman, {\it Units of Group Rings}, Proc. London Math. Soc. 46 (1971), 205--214.

%
%
%

\bibitem{Karpilovsky1} G. Karpilovsky, {\it Commutative Group Algebras}, Monographs and Textbooks in Pure and Applied Mathematics, Marcel Dekker Inc., New York, v. 78 (1983).

\bibitem{Karpilovsky18g} G. Karpilovsky, {\it Unit Groups of Classical Rings,} Oxford University Press, (1988).


\bibitem{Karpilovsky} G. Karpilovsky, {\it Units of Group Rings}, Longman Scientific \& Technical, (1989).

%
%
%
%

\bibitem{Passman} D. S. Passman, {\it  The Algebraic Structure of Group Rings}, John Wiley \& Sons, (1977).


%
%
\bibitem{Polcino} C. Polcino Milies and S. K. Sehgal, {\it An Introduction to Group Rings}, Kluwer Academic Pu\-blishers, (2002).
%
%
%
%
%


\bibitem{Sehgal0} S. K. Sehgal, {\it Topics in Group Rings}, Marcel Dekker, (1978).

\bibitem{Sehgal} S. K. Sehgal, {\it Units in Integral Group Rings},  Longman Scientific \& Technical, (1993).
%


\bibitem{Washington} L. C. Washington, {\it Introduction to Cyclotomic Fields}, Springer-Verlag (1997).


%
%
\bibitem[GAP]{GAP4}
The GAP~Group, {\it GAP - Groups, Algorithms, and Programming, Version $4.4.12$}; $2008,$ (http://www.gap-system.org).

\bibitem[MAPLE]{Maple}{\it Maple 13, MAC OS X Version}; $2009.$
%

\end{thebibliography}
\end{document}